\documentclass{amsart}
\usepackage{amsmath,amssymb, amsfonts}

\newcommand{\C}{\mathcal{C}}
\newcommand{\M}{\mathcal{M}}
\newcommand{\End}{\textnormal{End}}
\newcommand{\weight}{\textnormal{weight}}

\newcommand{\pr}{ \text{pr}}

 \newtheorem{theorem}{Theorem}

\theoremstyle{definition}
\newtheorem{example}[theorem]{Example}
\newtheorem{remark}[theorem]{Remark}

\title{Generalizations of the Goulden--Jackson Cluster Method}
\author{Elizabeth J. Kupin}
\address{SAS - Mathematics, Rutgers University, 110 Frelinghuysen Road, Piscataway, NJ  08854}
\email{ekupin@math.rutgers.edu}

\author{Debbie S. Yuster}
\address{DIMACS Center, Rutgers University, 96 Frelinghuysen Road, Piscataway, NJ 08854}
\email{yuster@math.rutgers.edu}
\thanks{Debbie Yuster was supported by a DIMACS Postdoctoral Fellowship.}

\begin{document}

\begin{abstract} We give several modifications of the Goulden--Jackson Cluster method for finding generating functions for words avoiding a given set of forbidden words. Our modifications include functions which can take into account various `weights' on words, including single letter probability distributions, double letter (i.e.\ pairwise) probability distributions, and triple letter probability distributions. We also describe an alternative, recursive approach to finding such generating functions. We describe Maple implementations of the various modifications. The accompanying Maple package is available at the website for this paper.
\end{abstract}
\maketitle
\section{Introduction}
Suppose we are given a finite alphabet, and a finite set of forbidden words in this alphabet. We would like to know how many $n$-letter words in our alphabet avoid the forbidden words as subwords or \emph{factors}, i.e.\ strings of consecutive letters. In order to do this we will find the generating function for the number of such words. In Section~\ref{naive} we describe a straightforward, recursive approach to solving this using ordinary generating functions. The remainder of this article deals with the Goulden--Jackson cluster method, a powerful method that considers overlap between forbidden factors in computing generating functions. The Goulden--Jackson cluster method was introduced in \cite{GJ1} and \cite{GJ2} and described very clearly and concisely in \cite{Noonan-Zeilberger}. For earlier work see~\cite{G-O}, and further extensions can be found in~\cite{Kong}. Applications of the Goulden--Jackson cluster method to genomics can be found in~\cite{Hao00}, \cite{Hao02}, and~\cite{Hao_language}. In Section~\ref{basic_GJ} we review the classical Goulden--Jackson cluster method. We then describe some modifications to the original Goulden--Jackson problem as follows:
 In Section~\ref{single} we modify the Goulden--Jackson cluster method to take single letter weights into account. In Section~\ref{double} we include double letter (i.e.\ pairwise) weights, and in Section~\ref{triple} we consider triple letter weights. Further generalizations are described in Section~\ref{others}.

All the methods discussed have been implemented in a Maple package which accompanies this paper. The Maple package, which includes documentation, can be found at the website for this paper
\footnote{The website for this paper is located at {\tt http://www.math.rutgers.edu/$\sim$ekupin/GJ.html}.}.

\section{The Straightforward Recursive Approach}\label{naive}

Given a finite alphabet $A$ and a set of forbidden or `bad' words $B$, we would like to find $a(n)$, the number of $n$-letter words in the alphabet $A$ that do not contain any members of $B$ as factors. Rather than find $a(n)$ directly, we will find the generating function \begin{equation}\label{genfcn}f(t)=\displaystyle\sum_{n=0}^{\infty}a(n)t^n   .\end{equation}

In order to find this generating function, we need not use the Goulden--Jackson cluster method, which will be described in Section~\ref{basic_GJ}. We can use a straightforward recursive approach, though as we will demonstrate, this method will not be as efficient as Goulden--Jackson. The approach contained in this section was described by Dr.\ Doron Zeilberger in his Spring 2008 Experimental Math class at Rutgers University.

We first illustrate the approach with an example. Suppose $A=\{a,b\}$ and $B=\{abb, ba\}$. We will start by decomposing the set of allowed words according to their first letter. Consider the set of allowed words beginning with $a$. Any such word is either $a$ itself, or consists of $a$ followed by a smaller word starting with either $a$ or $b$. What are the restrictions on the smaller word, following the initial letter $a$? If it starts with $a$, it must still avoid $abb$ and $ba$, but there are no additional restrictions. However, if the word following the initial letter $a$ begins with $b$, it must avoid $abb$ and $ba$, but in addition it must not begin with $bb$, so as to avoid the forbidden word $abb$. This gives a correspondence among different sets of words.

Denote by $[B, a_i, \{w_1,...w_n\}]$ the set of words avoiding members of $B$, starting with $a_i$, and avoiding any word in $\{w_1,...w_n\}$ as an initial subword. Translating the above example into this notation, we have:
\begin{equation}\label{naive1} [B,a,\{\}] \leftrightarrow \{a\} \cup [B, a,\{\}] \cup [B, b,\{bb\}],\end{equation} where the latter two terms on the right hand side describe the allowable subwords following the initial letter $a$.

Now let us consider the set of allowable words beginning with $b$. Either the word is $b$ itself, or consists of $b$ followed by a smaller word starting with either $a$ or $b$. Of course the letter following the initial $b$ cannot be $a$, since this would form the forbidden word $ba$, but rather than exclude this \emph{a priori}, we will instead say that any word following the initial $b$ and starting with $a$ must not begin with the word $a$. Of course no words satisfy this condition so the corresponding set will be empty.  Any word following the initial $b$ that starts with $b$ must avoid the forbidden words in $B$ but has no additional restrictions. To summarize, we have
\begin{equation}\label{naive2} [B,b,\{\}] \leftrightarrow \{b\} \cup [B,a,\{a\}] \cup [B,b,\{\}]. \end{equation}
We can express the allowable words with the decomposition
\begin{equation}\label{naive3} \left[B,*,\{\}\right] \leftrightarrow \{\text{empty\_word}\} \cup [B,a,\{\}] \cup [B,b,\{\}], \end{equation}
where $ \left[B,*,\{\}\right] $ denotes the set of all words in alphabet $A$ avoiding the words in $B$ (with no additional restrictions).

Eventually we will turn these correspondence relations into equations, by taking a weighted count of the set elements. First, in order to solve explicitly for the latter two sets in the above relation, we must decompose the remaining sets on the right hand sides of (\ref{naive1}) and (\ref{naive2}), as well as any new sets arising from those relations. This leads to the following:
\begin{eqnarray*}
\left[B, b,\{bb\}\right] &\leftrightarrow& \{b\} \cup [B,a,\{a\}] \cup [B,b,\{b\}]  \\
\left[B, a,\{a\}\right] &\leftrightarrow& \emptyset \\
\left[B, b,\{b\}\right] &\leftrightarrow& \emptyset .
\end{eqnarray*}
If one simply wants to know the number of allowable $n$-letter words, a generating function such as the one given in equation~(\ref{genfcn}) can be found. However, it is possible to find variant generating functions which give more information about the allowable words. These variants will be described in later sections. In order to find the various generating functions, we will make use of a weighted counting system, performing a weighted count of the words in the sets above.
Taking weights in (\ref{naive3}) gives the desired generating function. The particular `weight' used varies based on the method (to be described in the following sections) so we postpone further calculations. It is worth pointing out, however, that one must be careful not to merely sum the weights of the sets on the right hand sides of the correspondence relations. Rather, it is necessary to account for the weights of the truncated initial letters, as well as any transition weights that may arise. See Example~\ref{single_recursive_example} for further details.

The Maple code for this recursive method can be found in our accompanying Maple package under the function names {\tt RecursiveSingle}, {\tt RecursiveDouble}, and {\tt RecursiveProbDouble}. These functions implement the straightforward recursive analogues of the cluster method generalizations to be described in Sections \ref{single}, \ref{double}, and \ref{ProbDouble}, respectively.

\section{Basic Goulden--Jackson Cluster Method}\label{basic_GJ}
We borrow from \cite{Noonan-Zeilberger} in briefly reviewing the basic Goulden--Jackson cluster method, and encourage the reader to consult this source for a more detailed exposition.

In order to find the generating function given in equation~(\ref{genfcn}), we will do a weighted count of \emph{marked words}. A \emph{marked word} is a pair $(w;S)$, where $w$ is a word in the alphabet $A$, and $S$ is an arbitrary multiset whose entries are members of $Bad(w)$, the forbidden words of $B$ contained as factors in $w$. We allow repetition in $S$ since a word may contain several copies of a given forbidden word.
If no subset $S$ is specified, we assume it is the empty set. We define the weight of a marked word as $\weight(w;S)=(-1)^{|S|}t^{|w|}$, where $|S|$ is the cardinality of $S$ and $|w|$ is the length of $w$. The weight of a set of marked words is obtained by summing the weights of the marked words in the set.
The generating function from equation~(\ref{genfcn}) now becomes 

\begin{equation}\label{ft-1}f(t)=\displaystyle\sum_{w\in A^*}\sum_{S\subset Bad(w)}(-1)^{|S|}t^{|w|} ,\end{equation}
 where $A^*$ is the set of all words in the alphabet $A$.
 In order to see why this is valid, consider an arbitrary word $w$ in our alphabet. This word will appear as the first argument in $2^k$ marked words, where $k$ is the number of forbidden subwords contained in $w$. 
In equation~(\ref{ft-1}), we sum over all possible subsets of $Bad(w)$. Thus if $w$ contains no forbidden subwords, it will be counted exactly once in the above sum. If $w$ contains $k$ forbidden subwords, $k >0$, the number of times it will be counted in the above sum is:
$$\displaystyle\sum_{i=0}^k (-1)^i {k \choose i} = (1+(-1))^k=0.$$
Thus every allowable word is counted once, while words containing forbidden words are not counted. We have verified the equivalence of equation~(\ref{ft-1}) and our original generating function, given in equation~(\ref{genfcn}).

We call a marked word $(w;S)$ a \emph{cluster} if neighboring factors in $S$ overlap (i.e.\ are not disjoint) in $w$, and the forbidden words of $S$ span all of $w$.
For example, if our alphabet $A$ is $\{a,b,c\}$ and the set of forbidden words $B$ is $\{ba,aca\}$, then the marked word $(bacaca; \{ba,aca,aca\})$ is a cluster. The marked word $(bacaca; \{aca,aca\})$ is not a cluster because the factors in $S$ do not span all of $w$ (the first $b$ is not part of a factor in $S$), and the marked word $(acaba; \{aca,ba\})$ is not a cluster because the factors in $S$ do not overlap. We will denote the set of all (nonempty) clusters by $\C$.

We now decompose $\M$, the set of marked words, into three groups: the empty word, marked words beginning with a letter that is not part of any cluster, and marked words beginning with a cluster. We thus obtain the decomposition $$\M= \{\text{empty\_word}\} \space  \cup  A \M  \space \cup   \C \M,$$ where an element of $A\M$ consists of a single letter of the alphabet $A$ prepended to a marked word and an element of $\C\M$ consists of a cluster prepended to a marked word. Let $m$ be the number of letters in $A$. By taking weights on both sides of the preceding equation, we obtain
$$ \weight(\M) = 1 + mt\cdot \weight(\M) + \weight(\C)\weight(\M).$$ But $\weight(\M)$ equals $f(t)$, as shown in equation~(\ref{ft-1}), so by substituting and solving for $f(t)$ we get
\begin{equation}\label{ft-eqn}f(t)=\frac{1}{1-mt-\weight(\C)}\end{equation}
and it remains to solve for $\weight(\C)$, which we will call the \emph{cluster generating function}.

In order to find $\weight(\C)$, we partition the set of clusters $\C$ according to the first forbidden word of the cluster. Let $\C[v]$ denote the set of clusters starting with forbidden word $v$. Then $\C=\displaystyle\bigcup_{v\in B}\C[v]$, and  $\weight(\C)=\displaystyle\sum_{v\in B}\weight(\C[v])$.

In order to find $\weight(\C[v])$, we will further decompose $\C[v]$ as follows: consider a cluster in $\C[v]$. Either it consists of $v$ alone, or we can remove $v$ from the list of forbidden words in our marked cluster, and what remains will contain a smaller cluster, beginning with some bad word $u$ such that some initial subword of $u$ coincides with some final subword of $v$. For example, consider the cluster $(bacaca; \{ba,aca,aca\})$, where the alphabet and forbidden word set are as above. This cluster is in $\C[ba]$. Removing the initial forbidden word $ba$ leaves a new cluster $(acaca; \{aca,aca\})$ in $\C[aca]$. 
 Let $O(v,u)$ be the set of possible `overlaps' of $v$ and $u$, that is, all possible (nonempty) intersections of final subwords of $v$ with initial subwords of $u$. Each of these overlaps corresponds to a way that our cluster can have its first two forbidden words be $v$ and $u$, respectively.  To create the smaller cluster we will peel off exactly the part of $v$ that does not overlap with $u$. For any word $v$ and a final subword $r$ of $v$, $v\backslash r$ will denote the word obtained by chopping $r$ from the end of $v$. For example, $abcb\backslash cb=ab$.
 This leads to the decomposition
 $$\C[v] \leftrightarrow \{(v;\{v\})\} \cup \displaystyle\bigcup_{u\in B}\bigcup_{r\in O(v,u)}\Big(\C[u]\cdot(v \backslash r)\Big),$$
where $W_1\cdot W_2$ is the concatenation of $W_1$ with $W_2$.
 Taking weights, we obtain the following linear equations. Note that the sum over $u\in B$ is negative (i.e.\ multiplied by $-1$) in order to compensate for having reduced the number of bad words in our cluster by one (because we are calculating weights of clusters containing one fewer bad word than the clusters in $\C[v]$).

$$\weight(\C[v]) = \weight( (v;\{v\})) - \displaystyle\sum_{u\in B}\Big(\weight(\C[u])\cdot\sum_{r\in O(v,u)} \weight(v \backslash r)\Big)$$

We can explicitly calculate $O(v,u)$, and so by writing this equation for $\C[v]$ for all $v \in B$, we obtain a sparse system of $|B|$ linear equations in $|B|$ unknowns. Solving for the $\weight(\C[v])$ and summing them gives us $\weight(\C)$, which can then be substituted into equation~(\ref{ft-eqn}), giving the desired generating function.

\subsection*{Variations of the Basic Cluster Method}\label{variations}

By changing how the weight of a word is defined, we can alter the interpretation of the resulting generating function. In the following sections we present several such variations. All the variations keep track of how many words of each length avoid the set of forbidden words. The first variation, described in Section~\ref{single}, also takes into account how many times each letter appears in any given `good' word, by adding extra variables into the weight function. One possible use of this is to substitute probabilities for these variables, thus giving a generating function which takes into account a probability distribution on the alphabet. Several other variations mentioned later also take into account double letter, or pairwise weights. These are variables corresponding to each ordered pair of letters in the alphabet. This allows tracking of which consecutive letter combinations occur, and can also allow for double letter probabilities to be filled in. Similarly, Section~\ref{triple} has variables in the weight function corresponding to ordered triples of letters. In the sections that follow, we describe these modifications to the basic Goulden-Jackson cluster method in detail.

\section{Single Letter Weights}\label{single}

In order to keep track not only of how many words of a certain length avoid certain subwords, but also which letters these words contain, we will redo the basic cluster method, using a different weight enumerator. The variation discussed in this section was initially described in~\cite{Noonan-Zeilberger}. The weight of a marked word $(w;S)$ (where $w= w_1w_2\cdots w_k$) will be 
$$(-1)^{|S|} t^{|w|} x_{w_1} x_{w_2} \ldots x_{w_k}.$$ For example, $\weight(abcab; \{ \}) = t^5 (x_a)^2 (x_b)^2 x_c$.

As in the original application of the Goulden--Jackson method, we use the decomposition $\M = \{\text{empty\_word}\} \space \cup  A\M \cup \C\M $.  This leads to the following recursive formula for $\weight(\M)$: 
$$\weight(\M) = 1 + t\sum_{a \in A} x_a \weight(\M)+ \weight(\C)\weight(\M).$$ 
Note that, since we are keeping track of letter weights, the second term records not just how many letters are in $A$, but exactly which ones appear. Simplifying, we get
\begin{equation}\label{single1}\weight(\M)= \frac{1}{1 - t\displaystyle\sum_{a \in A} x_a - \weight(\C)}.\end{equation}

All that remains is to solve for the cluster generating functions, $\weight(\C)$. We do this exactly as in the original Goulden--Jackson cluster method, with the same decomposition:  $\weight(\C)= \sum_{v \in B} \weight(\C[v])$.
We will write the same system of linear equations as before, except that the weight function is different. In particular, we still have

$$\weight(\C[v]) = \weight( (v;\{v\})) - \displaystyle\sum_{u\in B}\Big(\weight(\C[u])\cdot\sum_{r\in O(v,u)} \weight(v \backslash r)\Big)$$
 for all $v \in B$, which becomes $$\weight(\C[v]) = -t^{|v|}x_{v_1}\ldots x_{v_{|v|}} - \displaystyle\sum_{u\in B}\Big(\weight(\C[u])\cdot\sum_{r\in O(v,u)} t^{|v\backslash r|} x_{v_1} \ldots x_{v_{|v \backslash r|}}\Big).$$

Solving for $\weight(\C[v])$ for each forbidden word $v$ and substituting back into equation~(\ref{single1}) yields the desired generating function.

\begin{example}\label{single_recursive_example}
Find the generating function of all words in the alphabet $\{a,b\}$ avoiding the forbidden words $abb$ and $ba$.\\
We will find the generating function in two ways: (1) using the cluster method described in this section, and (2) using the straightforward recursive approach from Section~\ref{naive}.

\begin{enumerate}
\item
We have:
\begin{eqnarray*}
\weight(\C[abb]) &=& -t^3x_ax_b^2 - \weight(\C[ba])t^2x_ax_b \\
\weight(\C[ba]) &=& -t^2x_bx_a - \weight(\C[abb])tx_b
\end{eqnarray*}
from which
\begin{eqnarray*}
\weight(\C[abb]) &=& \dfrac{-t^3x_ax_b^2 + t^4{x_a^2}{x_b^2}}{1-t^3{x_a}{x_b^2}}\\
\weight(\C[ba]) &=& -t^2x_bx_a + \dfrac{t^4{x_a}{x_b^3}-t^5{x_a^2}{x_b^3}}{1-t^3{x_a}{x_b^2}},
\end{eqnarray*}
and therefore $\weight(\C)= \dfrac{-t^2{x_a}{x_b} - t^3{x_a}{x_b^2} + t^4{x_a^2}{x_b^2} + t^4{x_a}{x_b^3}}{1-t^3{x_a}{x_b^2}}$. Substituting this into equation~(\ref{single1}) yields the desired generating function:
$$\weight(\M) = \dfrac{1-t^3{x_a}{x_b^2}}{1-t{x_a}-t{x_b}+t^2{x_a}{x_b}}.$$
Taking the first few terms of the Taylor expansion of this generating function yields:
$$1+({x_a}+{x_b})t+({x_a}{x_b}+{x_a}^2+{x_b}^2)t^2+({x_a}^2{x_b}+{x_a}^3+{x_b}^3)t^3+({x_a}^3{x_b}+{x_a}^4+{x_b}^4)t^4+O(t^5)$$
The constant term $1$ corresponds to the empty word. The coefficients of powers of $t$ correspond to the allowable words. For example, the coefficient of $t^3$ corresponds to the permissible 3-letter words $aab$, $aaa$, and $bbb$.

\item Returning to the notation of Section~\ref{naive}, we need to take the weight of $\left[B, *,\{\}\right].$ Recall the following set decompositions:
\begin{eqnarray*}
\left[B,*,\{\}\right] &\leftrightarrow& \{\text{empty\_word}\} \cup [B,a,\{\}] \cup [B,b,\{\}] \\
 \left[B,a,\{\}\right] &\leftrightarrow& \{a\} \cup \left[B, a,\{\}\right] \cup \left[B, b,\{bb\}\right] \\
 \left[B,b,\{\}\right] &\leftrightarrow& \{b\} \cup \left[B,a,\{a\}\right] \cup \left[B,b,\{\}\right] \\
\left[B, b,\{bb\}\right] &\leftrightarrow& \{b\} \cup [B,a,\{a\}] \cup [B,b,\{b\}]  \\
\left[B, a,\{a\}\right] &\leftrightarrow& \emptyset \\
\left[B, b,\{b\}\right] &\leftrightarrow& \emptyset .
\end{eqnarray*}
It remains to take weights of all the sets listed, from the bottom up, and solve for the unknown weights. We must be careful, however, to distinguish between identical sets on the left hand sides and right hand sides of the correspondence relations. For example, consider the correspondence $${\left[B,a,\{\}\right] \leftrightarrow \{a\} \cup \left[B, a,\{\}\right] \cup \left[B, b,\{bb\}\right]}.$$ The weight of the left hand side is simply $\weight(\left[B,a,\{\}\right])$, while the latter two sets on the right hand side are assumed to have had their initial letter $a$ removed. Thus, the total weight of the right hand side is $\weight(a) + \weight(a)\cdot\weight(\left[B, a,\{\}\right]) + \weight(a)\cdot\weight(\left[B, b,\{bb\}\right])$. Solving from the bottom up, we find:
\begin{eqnarray*}
\weight(\left[B, b,\{b\}\right]) &=& 0 \\
\weight(\left[B, a,\{a\}\right]) &=& 0 \\
\weight(\left[B, b,\{bb\}\right]) &=& \weight(b) \\
\weight(\left[B,b,\{\}\right]) &=& \weight(b) + \weight(b)\weight(\left[B,b,\{\}\right]) \\
\weight( \left[B,a,\{\}\right]) &=& \weight(a) + \weight(a)\weight(\left[B, a,\{\}\right])\\ &&+ \weight(a)\weight(\left[B, b,\{bb\}\right]) \\
\weight(\left[B,*,\{\}\right]) &=& \weight(\text{empty\_word})+\weight([B,a,\{\}]) + \weight([B,b,\{\}]) .
\end{eqnarray*}

Solving for each left hand side quantity and substituting into the last equation, which is the equation for $\weight(\M)$, we find:
\begin{eqnarray*}
\weight(\M) &=&1 + \dfrac{tx_a+t^2x_ax_b}{1-tx_a} + \dfrac{tx_b}{1-tx_b} \\
		&=& \dfrac{1-t^3x_ax_b^2}{1-t{x_a}-t{x_b}+t^2{x_a}{x_b}}.
\end{eqnarray*}
\end{enumerate}
\qed
\end{example}

We have implemented this modification of the original Goulden--Jackson cluster method, and the code is available in our accompanying Maple package under the function name {\tt SingleGJ}.

\section{Double Letter Weights}\label{double}

Sometimes we would like to keep track not just of how many times each letter appears in a word, but also which consecutive letter pairs appear. This could be relevant, for example, if studying English words, when the pair `QU' is many times more likely to appear than the pair `QB'. In order to keep track of such data, we introduce \emph{double letter weights}, that is, variables which represent the occurrence of consecutive letter pairs.

To include double letter weights, the weight of a marked word $(w;S),$ where $w=w_1w_2w_3 \cdots w_k$, will now be 

$$(-1)^{|S|} t^k (x_{w_1} \ldots x_{w_k})(x_{w_1, w_2} x_{w_2, w_3} \ldots x_{w_{k-1}, w_k}). $$ 
We will denote this new weight function $W((w;S))$. For example,  $W((cat; \{\})) = t^3 (x_c x_a x_t)(x_{c,a} x_{a,t})$. This new weight function does not have all of the nice properties of weight functions we have seen in the earlier methods. In particular, concatenation of words no longer corresponds to a simple multiplication of weights. To see why this is true, consider the word $abab$ as the result of concatenating $ab$ with itself. In this case, $W((uu;\{\}))$ does not equal $W((u;\{\}))^2$. $W((ab; \{\}))=t^2(x_a x_b)(x_{a,b})$ and so $W((ab;\{\}))^2=t^4(x_a)^2 (x_b)^2(x_{a,b})^2,$ while 
$W((abab;\{\})) = t^4(x_a)^2 (x_b)^2(x_{a,b})^2 x_{b,a}.$

In general, whenever we concatenate two strings we need to account for the double letter weight that crosses from one string to the next. We call this the extra factor the \emph{transition weight}. The original cluster method involves decomposing $\M$, then using the fact that a disjoint union of sets corresponds to addition of weight functions, and concatenation corresponds to multiplication. We can still use this basic principle, but we must be more careful with concatenation. In particular, whenever we concatenate strings we will need to know the last letter of the first string and the first letter of the second string, in order to be able to multiply by the appropriate transition weight. This forces us to change how $\M$ is decomposed.

In the original method, we used the decomposition  $$\M= \{\text{empty\_word}\} \space  \cup  A \M  \space \cup   \C \M.$$ This involves concatenation in two places: in the second term we concatenate an arbitrary marked word to a single letter, and in the third term we concatenate an arbitrary marked word to a cluster. To incorporate the transition weights we will need to know the first letter of an arbitrary marked word, as well as the last letter of an arbitrary cluster. 

We start by splitting up the set of marked words according to their first letter. Let $\M_a$ be the set of marked words that start with $a$, $\M_b$ be the set of marked words that start with $b$, and so on. We have $$\M = \{ \text{empty\_word} \} \cup \left( \bigcup _{a \in A} \M_a \right).$$ 

To find $\weight(\M_a)$, we examine the different types of marked words that can begin with the letter $a$. Such a word may be $a$ itself, or we can peel off the initial $a$ to get a shorter marked word (assuming the initial $a$ is not part of a cluster), or the word begins with a cluster that begins with $a$. Let $B_a$ be the set of forbidden words beginning with $a$. We have the decomposition 
$$\M_a = {a} \cup \left( \bigcup _{b \in A} a\M_b \right) \cup \left( \bigcup _{v \in B_a} \C[v]\M \right). $$
Accounting for the fact that the entire marked word may be a cluster, we get
\begin{equation}\label{Ma decomp.}\M_a = {a} \cup \left( \bigcup _{b \in A} a\M_b \right) \cup \left( \bigcup _{v \in B_a}  \bigcup _{b \in A} \C[v] \M_b \right) \cup \left( \bigcup _{v \in B_a} \C[v] \right).\end{equation}
In this manner we keep track of the first letter of each marked word. It remains to address concatenation in the cluster generating functions.

\subsection* {Cluster Generating Functions}

The decomposition of $\C[v]$ in the basic cluster method is based on the idea that if we have a cluster beginning with a bad word $v$, the cluster is either just that word, or we can peel the first word off and get a smaller cluster beginning with a bad word $u$ that has some non-trivial overlap with $v$. Thus we have
$$\C[v]= v \cup \left(  \bigcup_{u \in B} \bigcup_{r \in O(v,u)} (v\backslash r)\C[u]  \right).$$

Since we are computing this for a specific  $v$, we know what the last letter of $v\backslash r$ will be. Moreover, we know what the first letter of $u$ will be, so the transition weight is easy to write down. Taking weights on both sides, we get
$$W(\C[v]) = W((v;\{v\})) - \left(\sum_{u\in B}\sum_{r \in O(v,u)} W(v \backslash r)\underbrace{x_{v_{|v\backslash r|},v_{|v\backslash r|+1}}}_{\text{transition weight}}W(\C[u]) \right).$$ 

The cluster method is based on computing weights of individual letters and clusters, then computing weights of marked words in terms of the letters and clusters they contain. However, computing weights of such concatenations requires the inclusion of transition weights, since we are incorporating double letter weights into our weight-enumerators. In order to do this, it becomes necessary to keep track of the last letter of each cluster. When computing cluster weights, we successively remove leading forbidden words from a cluster, until we are left with a cluster consisting of only one forbidden word. By keeping track of its last letter, we are keeping track of the last letter of the original cluster. Thus, it suffices to record the last letter of single-word clusters only. We do this by adding a dummy variable to the end of each one-word cluster. This dummy variable records the last letter of a one-word cluster:

$$W(\C[v]): = W((v;\{v\}))\underbrace{\End_{v_k}}_{\text {dummy}}-  \left(\sum_{u\in B}\sum_{r \in O(v,u)} W(v \backslash r) x_{v_{|v\backslash r|},v_{|v\backslash r|+1}}W(\C[u]) \right).$$ 

Now we have a system of linear equations with variables $W(\C[v])$, for $v \in B$, and we can solve for each of these in terms of the dummy variables $\End_a$, where $a \in A$. However, we don't want our final equation in terms of these variables.
When we prepend a cluster to an arbitrary word beginning with $a$, and wish to take the resulting weight, we must first replace all occurrences of $\End_b$ (for any letter $b$) with the transition weight $x_{b,a}$.
 This gives us a way to calculate the transition weight directly from the cluster generating function, allowing us to use equation~(\ref{Ma decomp.}). We cannot write down exactly what the transition weight will be in general, since it will depend on which cluster we are looking at, but if we let $T_{\C}$ denote the transition weight calculated for a specific cluster, we have

$$W(\M_a) = tx_a + \sum_{b\in A} t x_a x_{a,b} W(\M_b) + \left(\sum_{f \in B_a}\sum_{c \in A}W(\C[f])\,T_{\C}\,W(\M_c)\right) + \sum_{f \in B_a} W(\C[f]).$$ 
 
 Now that we can calculate the $W(\M_a)$, we can put them together to find $W(\M)$.
The Maple code for the cluster method, implementing the weight enumerator described in this section, can be found in our accompanying Maple package under the function name {\tt DoubleGJ}.

\begin{example} We return to the setup from Example~\ref{single_recursive_example}, redoing the example with double letter weights. Recall the problem: find the generating function of all words in the alphabet $\{a,b\}$ avoiding the forbidden words $abb$ and $ba$.

Running {\tt DoubleGJ} to get the generating function and doing a Taylor expansion, we see the first several terms are:
\begin{multline*}1+(x_b+x_a)t+(x_ax_{a,b}x_b+x_a^2x_{a,a}+x_{b,b}x_b^2)t^2+(x_a^2x_{a,a}x_{a,b}x_b+x_a^3x_{a,a}^2+x_{b,b}^2x_b^3)t^3\\
+(x_a^3x_{a,a}^2x_{a,b}x_b+x_a^4x_{a,a}^3+x_{b,b}^3x_b^4)t^4+(x_a^4x_{a,a}^3x_{a,b}x_b+x_a^5x_{a,a}^4+x_{b,b}^4x_b^5)t^5+O(t^6).
\end{multline*}
The coefficient of $t^4$ has terms corresponding to the allowable four-letter words $aaab$, $aaaa$, and $bbbb$.
\qed

\end{example}

\begin{remark}[Comparison of straightforward recursive approach and Goulden--Jackson cluster method]
The straightforward recursive approach will usually require a system of at least $|A|+|B|$ equations and unknowns, often more (where $A$ is the alphabet and $B$ is the set of forbidden words). The Goulden--Jackson method with double letter weights first solves a system of size $|B|$ (the cluster generating functions), and then a system of size $|A|$. Even if the number of equations is roughly the same, by breaking things apart a bit we would still expect the Goulden--Jackson method to be slightly faster. In practice, however, the differences seem to be small for small examples.
In general, the best approach depends on the situation.
 If there are many forbidden words with a lot of overlap, then Goulden--Jackson may take longer to compute the cluster generating functions, making it slower. However, if there are fewer, but longer, forbidden words, the straightforward recursive approach may require many more than $|A|+|B|$ equations and therefore take longer.

\end{remark}

\section{Triple Letter Weights}\label{triple}

As a further generalization of the original cluster method, in this section we will keep track of the occurrences of each letter in a word, the occurrences of consecutive letter pairs, and the occurrences of consecutive letter triples.
We will use a new weight function, $W',$ that accounts for all single letters, letter pairs and letter triples in a word. For example, $$W'((abcabc; \{\}))=t^6(x_a^2 x_b^2 x_c^2)(x_{a,b}^2 x_{b, c}^2 x_{c,a})(x_{a,b,c}^2 x_{b,c,a} x_{c,a,b}).$$

Note that, as in Section~\ref{double}, we do not have $W'((uv;\{\}))={W'((u;\{\}))\cdot W'((v;\{\}))}$. In fact in the example above we can see that $W'((abcabc; \{\}))$ has three extra terms that do not appear in $W'((abc;\{\}))^2$. The term $x_{c,a}$ is from the double letter transition weight, as described in Section~\ref{double}. The other two extra terms, $x_{b,c,a}$ and $x_{c,a,b}$, correspond to the triples that cross between the two factors. These are also considered transition weights, and our main problem in this section will be modifying our methods so that it is possible to figure out exactly what these transition weights will be. 
 
For considering double letter weights but not triple letter weights, it is necessary to know
the last letter of the first string and the first letter of the second string when concatenating, in order to write down the appropriate transition weight. Now that we are also considering letter triples, we need to know the last two letters of the first string and the first two letters of the last string in order to get both triple transition weights. For example, to concatenate $abcd$ and $klmn$, the transition weights will come from the strings $cdk$, $dkl$, and $dk$, so in order to find these weights we need the final two letters of $abcd$, as well as the first two letters of $klmn$. 
 
Recall that in our original setup, we decomposed the set of marked words $\M$ as follows: $$\M= \{\text{empty\_word}\} \space  \cup  A \M  \space \cup   \C \M.$$ In this decomposition, we append marked words to individual letters, and marked words to clusters. In order to include triple letter weights we will need to know the first two letters of an arbitrary marked word.
We will also need the last two letters of an arbitrary cluster.
 
In order to keep track of the first two letters of an arbitrary marked word, let $\M_{ab}$ be the set of marked words beginning $ab$, and decompose $\M$ as follows: $$\M = \{\text{empty\_word}\} \cup S \cup \left( \bigcup_{a,b \in A} \M_{ab} \right).$$ Here $S$ is the set of one-letter marked words. Taking weights gives us 

\begin{equation}\label{W} W'(\M)= 1+ \sum_{a \in A} W'((a;\{\}))+ \sum_{a \in A}\sum_{b \in A} W'(\M_{ab}).\end{equation}

To solve for $W'(\M_{ab}),$ we decompose $\M_{ab}$ further. A marked word in $M_{ab}$ could be the two-letter marked word $(ab; \{\})$, or the single letter $a$ followed by an arbitrary marked word beginning with $b$, or it could consist of a cluster beginning with $ab$, followed by an arbitrary marked word. Let $B_{ab} \subset B$ be the set of all forbidden words that begin with $ab$. We may assume that there are no one-letter forbidden words (in that case we would simply remove that letter from the alphabet). Thus B is completely partitioned into the $B_{ab}.$ We have

$$\M_{ab} = \{ab\} \cup \left(\bigcup_{c \in A} a\M_{bc}\right) \cup \left( \bigcup_{v \in B_{ab}} \C[v]\M \right).$$
  By taking weights on both sides, using the substitution given in equation~(\ref {W}), and adding in the appropriate transition weights when we can, we get 
  
  \begin{multline}W'(\M_{ab})=W'(ab) + \sum_{c \in a} W'(a) \underbrace{x_{a,b} x_{a,b,c}}_{\text{transition}} W'(\M_{bc}) + \sum_{v \in B_{ab}} W'(\C[v])\\+\sum_{v \in B_{ab}} \sum_{c \in A} W'(\C[v])\,T_{\C}\, W'(c)+ \sum_{v \in B_{ab}} \sum_{c,d \in A} W'(\C[v])\,T_{\C}\, W'(\M_{cd}) \end{multline}

We aren't yet able to fill in the transition weights $T_{\C}$ when they occur at the end of a cluster. This is the same problem that occurs in the double letter case, and the solution is the same. Into the cluster generating function we will put the dummy variables $\End_{ab},$ for all $a,b \in A$. They will go in the same place $\End_a$ went, and record the last two letters of a cluster. Then we can solve for the $W'(\M_{ab})$ in terms of the dummy variables and substitute for them when needed.

We have implemented the cluster method incorporating triple letter weights in the function {\tt TripleGJ}. The code can be found in our accompanying Maple package.

\begin{remark}
While it is possible to further generalize this method to include weights for four-letter strings or higher, 
we can see even from triple letter weights why this might be problematic. For one thing, efficiency suffers greatly. Using only double letter weights leads to a system of $|A|$ linear equations, while including triple letter weights requires $|A|^2$ equations. In general, including weights of strings up to length $k$ leads to a system of $|A|^{k-1}$ linear equations. 

Moreover, small forbidden words become increasingly problematic. In the triple letter case, we relied on the fact that having a one-letter forbidden word is equivalent to looking at a smaller alphabet. If we wanted to weight four-letter strings, we would partition the forbidden words based on their first three letters, so two-letter forbidden words would be a problem. As we include longer and longer subwords in our weight function, we get more and more short forbidden word exceptions.
\end{remark}

\section{Further Generalizations}\label{others}

In this section we introduce further variations of the original Goulden--Jackson cluster method. These variations have been implemented in our accompanying Maple package, and we refer to the variations according to their corresponding function names in our software.

\subsection{ProbDoubleGJ, ProbTripleGJ}\label{ProbDouble}

{\tt ProbDoubleGJ} and {\tt ProbTripleGJ} are variants of {\tt DoubleGJ} and {\tt TripleGJ} (respectively), that are specifically designed for applications with an underlying Markov chain structure. We can think of the set of states in a Markov chain as an alphabet, and a word in that alphabet will correspond to a history of steps in the Markov process. We move to the next state (or equivalently, add the next letter to our word) with some probability that depends only on the current state.

{\tt ProbDoubleGJ} returns a generating function for subword avoidance in which each word is weighted with all of the letter pairs it contains, as well as the single letter weight for the initial single letter only. For example,  $\weight((cat; \{\})) = t^3 (x_c)(x_{c,a} x_{a,t})$. We can interpret the double letter weight as a conditional probability: $x_{a,b}$ is the probability we will move to state $b$, given that we are currently at state $a$. The initial single letter probability represents the probability of starting in a given state. 

{\tt ProbTripleGJ} yields a generating function for subword avoidance where each word is weighted according to the triple letter strings it contains, as well as its initial double letter pair. The weight of a word of length one is the single letter weight. For example, $$\weight((abcabc; \{\}))=t^6(x_{a,b})(x_{a,b,c}^2 x_{b,c,a} x_{c,a,b}).$$ We can interpret the triple letter weight $x_{a,b,c}$ as the conditional probability of seeing $c$ given that the letters immediately preceding it are $ab$, and the initial double letter weight as the probability of satisfying a two-state initial condition. 

Note that these programs aren't contained in the original double letter weight and triple letter weight programs, in the sense that we can't get all the results here by a clever choice of the weights in those programs (described in Sections~\ref{double} and \ref{triple}). Our goal here is to selectively weight only the single or double letters that show up at the beginning of a word, whereas in the other programs, we assigned a weight to every letter in a word, regardless of where it appeared. Nevertheless, with some simple modifications to the original programs, we can get the desired results.

In order to modify our original double letter weight method (implemented in the function {\tt DoubleGJ}), we need only change the weight enumerator used. The new weight enumerator will consist of the initial single letter weight multiplied by $V(w;S)$, which we define to be the product of all consecutive letter pair weights. Define $V(a;\{\})=t$ for all single letter words $a$. For example, $V(abba;\{\})= x_{a,b}x_{b,b}x_{b,a}.$
We first solve for the $V(\M_a)$ as in Section~\ref{double}, yielding generating functions that enumerate all words beginning with $a$, and weighted only by their double letter components. It remains to multiply by initial letter weights and put things back together:
$$f(t)=1+\sum_{a \in A} x_aV(\M_a).$$

The situation is similar in adapting the triple letter weight method. We define $V'(w,S)$, which counts triple letter weights, i.e.\ $V'(abab;\{\})=x_{a,b,a}x_{b,a,b},$ $V'(ab;\{\})=t^2,$ $V'(a;\{\})=t$. By solving for the $V'(\M_ab)$ as in Section~\ref{triple}, we get generating functions for marked words, weighted only by their triple letter weights. Multiplying by the initial double letter weight and putting everything together yields:
$$f(t)=1+\sum_{a \in A} x_at+ \sum_{a,b \in A}x_{a,b}V'(\M_{ab}).$$

\begin{example}
Suppose we are given all initial letter probabilities and pairwise transition probabilities. For example, let's start with the alphabet $\{a,b\}$, initial letter probabilities $x_a=0.75$, $x_b=0.25$, and double letter probabilities $x_{a,a}=0.5$, $x_{a,b}=0.5$, $x_{b,a}=0.7$, $x_{b,b}=0.3$.

Suppose we would like to find the probability of avoiding the forbidden words $bbb$ and $ab$.
Running {\tt ProbDoubleGJ} produces the generating function
$$f(t)=-\dfrac{1}{100}\cdot\dfrac{100t+200+3t^3+25t^2}{t-2}.$$
The first few terms of the Taylor expansion are
$$1+t+\frac{5}{8}t^2+\frac{131}{400}t^3+\frac{131}{800}t^4+O(t^5).$$
The coefficient of $t^2$ is the probability of seeing one of the three allowable two letter words ($aa$, $ba$, or $bb$). We can calculate this probability by subtracting the probability of $ab$ from $1$. The probability of $ab$ is the initial probability of $a$, $0.75$, multiplied by the digraph probability of $ab$, $0.5$. Given a two letter word, we see that the probability of seeing the forbidden word $ab$ is 3/8, thus the probability of an allowable two letter word is $5/8$.
\qed

\end{example} 

\begin{example}[Modeling the English language]\label{ex:monkey}

We will consider a passage of written English as a long string over an alphabet of 27 characters: the 26 lowercase letters of the alphabet, and a character $SP$ that corresponds to a space. To keep things simple, we ignore punctuation and capitalization. In this example we will use the words  \emph{character} and \emph{string} when referring to the model and including the character $SP$, and \emph{letter} and \emph{word} when referring to English.

To use \texttt{ProbDoubleGJ}, we need to know information about the frequencies of single characters and character pairs. Where this information comes from can have a huge effect on the accuracy of this model. For details on how we obtained the letter frequencies used in this example, we refer the reader to Appendix~\ref{app:monkey}.

In addition to the table of frequencies, we will need a set of forbidden strings. For this example, we will use the forbidden string ``$SP,t,h,e$''. This corresponds to typing a word beginning with ``$the$'', a common word beginning. 

Running \texttt{ProbDoubleGJ} outputs a rational function with degree 29 polynomials in both the numerator and the denominator. We will call this function $F(x)$, and we list the first several terms of its Taylor expansion: 
\begin{multline*}F(x)=1+x+x^2+x^3+0.9992162308\cdot x^4+0.9992162308\cdot x^5 +0.9991288963\cdot x^6\\
 +0.9990341113\cdot x^7 +0.9989403663\cdot x^8+0.9988466147\cdot x^9+0.9987529643\cdot x^{10}\\
+0.9986592740\cdot x^{11} +0.9985656098\cdot x^{12}+0.9984719485\cdot x^{13}+0.9983782981\cdot x^{14}\\
 +0.9982846557\cdot x^{15}+0.9981910224\cdot x^{16} +0.9980973977\cdot x^{17}+0.9980037819\cdot x^{18}\\
 +0.9979101748\cdot x^{19} +0.9978165765\cdot x^{20}+0.9977229870\cdot x^{21} +0.9976294063\cdot x^{22}\\
 +0.9975358344\cdot x^{23}+0.9974422712\cdot x^{24} +0.9973487168\cdot x^{25}+ O(x^{26})
\end{multline*}The coefficient of $x^n$ corresponds to the probability of avoiding ``$SP,t,h,e$'' in a length $n$ string, according to the probability distributions given. It makes sense that the coefficients of $x$, $x^2$, and $x^3$ are all one, as the forbidden string is four characters long. Looking forward in the series, the coefficient of $x^{100}$ is $0.9903570875$, the coefficient of $x^{200}$ is $0.9811110978$, the coefficient of $x^{300}$ is $0.9719514288$, and the coefficient of $x^{400}$ is $0.9628772746$.

This means that if a monkey is banging keys on a typewriter according to the probability distributions given, even after 100 keystrokes (including the spacebar) the monkey only has a 1\% chance of typing a word beginning ``the...''. This probability climbs to just under 4\% after typing $400$ keystrokes.

To compare, and to answer the age old question about monkeys being able to type Hamlet, according to our model, the probability that a monkey could type ``to be or not to be'' after $300$ keystrokes is $5.861724357\cdot10^{-21}$. In comparison, if the monkey were hitting keys at random, with a $1/27$ probability of typing each letter or the space key, then the probability of typing ``to be or not to be'' in the first $300$ keystrokes is $6.62874079\cdot10^{-27}$.
More examples are available on the website for this paper, as well as some of the functions used to analyze these long strings of data, and the dictionary used in this example.
\qed
\end{example}

\subsection{DoubleGJIF}

The program {\tt DoubleGJIF} is a generalization of {\tt ProbDoubleGJ}. This function returns the generating function for subword avoidance where each word is weighted by all of its digraph weights as well as the initial single letter and final single letter weights (the {\tt IF} stands for Initial-Final). The implementation allows for setting different values for the weight of a single letter depending on whether it is the first or the last letter in a word. If we set all of the final letter probabilities to $1$, the program reduces to {\tt ProbDoubleGJ}.

There is a way to modify {\tt ProbDoubleGJ} to obtain {\tt DoubleGJIF} - simply multiply by the final letter weights as they occur. Of course, locating the end of a word in the decomposition is a bit more complicated, but it turns out that there are only a few places where we need to add terms.

Recall that in {\tt ProbDoubleGJ} we have the equation $$f(t)=1+\sum_{a \in A}x_{a}V(\M_{a}),$$  and adapted from {\tt DoubleGJ} we have $$V(\M_a) = t + \sum_{b\in A} t x_{a,b} V(\M_b) + \left(\sum_{f \in B_a}\sum_{c \in A}V(\C[f])\,T_{\C}\,V(\M_c)\right) + \sum_{f \in B_a} V(\C[f]),$$ where $V$ is the weight enumerator that weights words only with the letter pairs they contain.  

Suppose a marked word ends with a single letter (as opposed to a cluster). We decompose the word by peeling off letters or clusters from the beginning of the word. Once we arrive at the final letter it will seem as if we are looking at a marked word consisting of a single letter. Therefore, if we multiply the term in the equation above that corresponds to a single letter with the appropriate final letter weight, we will have successfully modified the end of every marked word that ends in a single letter.

Similarly, in dealing with the marked words that end in a cluster, we need only change the term in the above sum that corresponds to a marked word made of one single cluster. In order to implement {\tt DoubleGJ} we had to locate the end of the clusters, and so they are now flagged with the dummy variables $\End_a, a \in A$. Normally, when we concatenate the empty word to a cluster (effectively ending a word with a cluster), we substitute $1$ for the terms $\End_a$, for all $a \in A$. Instead, we can substitute the final letter probabilities for $\End_a$, and we will have successfully modified all the marked words that end in a cluster.  Since every nonempty marked word must end with a single letter or a cluster, we have added in the final letter probability to every marked word.

\subsection{DoubleGJst}

All of the programs we have discussed so far return a generating function in terms of $t$, but in fact any of the programs can be modified so that they return a generating function in two variables, $s$ and $t$:
$$f(s,t)=\sum_{n=0}^{\infty} \sum_{i=0}^n a(i,n)s^i t^n,$$ where $a(i,n)$ is the number (or weight, depending on the application) of words of length n that contain exactly $i$ forbidden subwords (counted with multiplicity). For example, if $A=\{a,b\}$ and $B=\{abb,ba\}$, then three out of the sixteen four-letter words contain no subwords in $B$: $aaaa$, $aaab$, and $bbbb$. Ten of the words contain exactly one forbidden subword: $aaba$, $aabb$, $abaa$, $abab$, $abbb$, $baab$, $baaa$, $bbaa$, $bbab$, and $bbba$. Finally, three of the words contain exactly two forbidden subwords each: $baba$, $babb$, $abba$. Therefore, we have $a(0,4)=3,$ $a(1,4)=10,$ and $a(2,4)=3$. Since this accounts for all 16 of the subwords of length 4, we must have that $a(3,4)=a(4,4)=0$. Therefore, the coefficient of $t^4$ in $f(s,t)$ will be $3+10s+3s^2$. In more complicated programs, $a(i,n)$ would give the weight of these words, not just the number of them. Not surprisingly, we can get this generating function by modifying the weight enumerator in any of the above programs.

Suppose we have a word $w$ that contains $k$ forbidden subwords. The word $w$ will be the first argument in $2^k$ marked words -- one for each subset of the $k$ forbidden subwords. In the original Goulden--Jackson method, we multiplied the weight of the letter in $w$ by $(-1)^{|S|}$, so that the number of times a word with $k$ forbidden subwords will be be counted is
$$\displaystyle\sum_{i=0}^k (-1)^i {k \choose i} = (1+(-1))^k=0^k,$$ 
which is $0$ if $k>0$ and 1 if $k=0$. If we would like to keep track of the number of forbidden subwords a word contains, we can simply replace the $(-1)$ in the weight function with $(s-1)$.
Thus the weight of a marked word $(w;S),$ where $w=w_1w_2w_3 \cdots w_k$, will be 
$$(s-1)^{|S|} t^k (x_{w_1} \ldots x_{w_k})(x_{w_1, w_2} x_{w_2, w_3} \ldots x_{w_{k-1}, w_k}). $$ 
 Under this new weight function, the number of times a word containing $k$ subwords will be counted is
$$\displaystyle\sum_{i=0}^k (s-1)^i {k \choose i} = (1+(s-1))^k=s^k,$$ 
as desired. The notion of including an extra variable to count the number of forbidden word occurrences was introduced in~\cite{Noonan-Zeilberger}. The functions {\tt SingleGJst}, {\tt DoubleGJst}, and {\tt ProbDoubleGJst} in our accompanying Maple package are the respective analogues of {\tt SingleGJ}, {\tt DoubleGJ}, and {\tt ProbDoubleGJ} incorporating the variable $s$.

\vspace{0.5cm}

\noindent {\bf Acknowledgements:}
This work began as a final project in Dr.~Doron Zeilberger's Experimental Mathematics course, given at Rutgers University in Spring 2008. The authors wish to thank Dr.~Zeilberger for suggesting this project and for his helpful advice and guidance.

\newpage
\appendix
\section{Data for Example~\ref{ex:monkey}}
\label{app:monkey}
 In order to obtain frequencies for single letter occurrences as well as pairwise letter frequencies, we analyzed a list of $20,422$ distinct English words that we will refer to as the dictionary. We created a large transition probability matrix by taking the frequency of a character pair and then dividing by the number of occurrences of the initial character. For example, the following table corresponds to the probabilities that any of the 27 characters should follow $a$:

\begin{tabular}{lll}
\pr(a,a)=0.00037487,& \pr(a,b)= 0.044235,& \pr(a,c)= 0.059454,\\
\pr(a,d)= 0.042885, & \pr(a,e)= 0.0030739,& \pr(a,f)= 0.010046 \\
\pr(a,g)= 0.033138,& \pr(a,h)= 0.0039736, &\pr(a,i)= 0.029090\\
\pr(a,j)= 0.00067476,& \pr(a,k)= 0.011771, & \pr(a,l)= 0.11553\\
\pr(a,m)= 0.039061,& \pr(a,n)=0.14342,& \pr(a,o)= 0.00074974\\
\pr(a,p)= 0.034638,& \pr(a,q)= 0.00089969,& \pr(a,r)=0.12003\\
\pr(a,s)= 0.052856, & \pr(a,t)=0.14530, & \pr(a,u)= 0.019793\\
\pr(a,v)= 0.014020,& \pr(a,w)= 0.0099715,& \pr(a,x)= 0.0044235\\
\pr(a,y)= 0.015070,& \pr(a,z)= 0.0044235,& \pr(a,SP)= 0.041086\\
\end{tabular}

The last value, $\pr(a, SP)$, was computed not from a letter pair in English, but from the frequency of words in English that end with $a$. The overall sum is $1$, since every occurrence of the letter $a$ is either followed by another letter, or occurs at the end of a word. 

Relatively speaking, the probabilities for the character $a$ are fairly evenly distributed. To compare, the corresponding values for $q$ look quite different:

\begin{tabular}{lll}
\pr(q,a)=0.0,& \pr(q,b)= 0.0,& \pr(q,c)= 0.0\\
\pr(q,d)= 0.0, & \pr(q,e)= 0.0,& \pr(q,f)= 0.0 \\
\pr(q,g)= 0.0,& \pr(q,h)= 0.0, &\pr(q,i)= 0.0\\
\pr(q,j)= 0.0,& \pr(q,k)= 0.0, & \pr(q,l)= 0.0\\
\pr(q,m)= 0.0,& \pr(q,n)=0.0,& \pr(q,o)= 0.0\\
\pr(q,p)= 0.0,& \pr(q,q)= 0.0,& \pr(q,r)=0.0\\
\pr(q,s)= 0.0, & \pr(q,t)=0.0, & \pr(q,u)= 0.99708\\
\pr(q,v)= 0.0,& \pr(q,w)= 0.0,& \pr(q,x)= 0.0\\
\pr(q,y)= 0.0,& \pr(q,z)= 0.0,& \pr(q,SP)= 0.0029240\\
\end{tabular}

Of the 342 occurrences of the letter $q$ in our dictionary, 341 of them are followed by the letter $u$, and exactly one of them is at the end of the word. Scrabble fanatics will no doubt appreciate that our dictionary is incomplete. 
We further remark that our model ignores context, and the fact that some words are more common than others in written English.

The last row of the probability matrix will be the probalities $\pr(SP,a)$, $\pr(SP, b)$, etc. These are the initial letter probabilities, in other words, the probability a word begins with a particular letter.
As a default, we set $\pr(SP,SP)=0$, ensuring that between two words there will only be one space. Calculating the single character frequencies is straightforward for the characters that are letters. We set the frequency of $SP$ according to the number of words in the dictionary, with the idea that between each word there must be exactly one space.

\bibliography{GJbib}
\bibliographystyle{alpha}

\end{document}